# Evaluating the Infinite

Toby Ord*
University of Oxford

I present a novel mathematical technique for dealing with the infinities arising from divergent sums and integrals. It assigns them fine-grained infinite values from the set of hyperreal numbers in a manner that refines the standard theories of summation and integration. This has implications in statistics (helping us work with distributions whose mean or variance is infinite), decision theory (allowing comparison of options with infinite expected values), economics (allowing evaluation of infinitely long streams of utility without discounting), and ethics (allowing evaluation of infinite worlds). There are even implications for finite cases, as the ability to handle these infinities undermines a common argument for bounded utility and the discounting of future utility.

Decision theory, economics, and ethics all encounter serious technical problems when the value of something goes to infinity. This could be a prospect (such as the St Petersburg gamble), an unending stream of benefits, or a population. In all cases, we want to evaluate or compare things with infinitely many finite parts, even when the sum or integral diverges to infinity. These technical problems for infinite cases have led many in these fields to adopt controversial workarounds (such as bounded utility or pure time preference) which change how we rank alternatives even in finite cases. This is premature, for there are technical solutions to the technical problems.

Standard theories of summation and integration use the extended reals, which aren't rich enough represent our intuitions. They have only one positive infinite value ($+\infty$) which forces all infinitely valuable things to be equal. The ordinal and cardinal numbers are also inadequate: as generalisations of natural numbers they are missing too much of the structure of the reals. But there is a very natural extension of the reals to infinite values (the *hyperreal numbers*) which is rich enough to allow a nonstandard theory of summation and integration that assigns fine-grained infinite values to the relevant sums and integrals.

## The hyperreal numbers

The hyperreals ($^*\mathbb{R}$) are an extension of the real numbers that includes both infinitesimal and infinite numbers (Robinson 1966, Keisler 1976, Goldblatt 1998). Just as real numbers can be defined in terms of sequences of rational numbers, so the hyperreal numbers can be defined in terms of sequences of real numbers. But there are two differences. We drop the requirement that the sequences are Cauchy-convergent (divergent sequences will allow for infinite numbers), and we use a finer grained distinction for when two sequences represent the same number (allowing for numbers that are infinitesimally different from each other).

---

* I would like to thank Marcus Pivato, Owen Cotton-Barratt, Kenny Easwaran, Alan Hájek, Jeffrey Russell, Fin Moorhouse, and Matthew van der Merwe for their very helpful discussions and comments. I'm especially grateful to Marcus, who has done much to pioneer this area and whose advice on this project over many years now has been invaluable.

For sequences to represent the same number, they don't merely have to approach each other. They have to be exactly equal at a 'large' set of indices. The set of all indices ($\mathbb{N}$) is considered large, as are all subsets of $\mathbb{N}$ that are only missing finitely many elements. So if two sequences are equal from some point onwards, they will be equal at a large set of indices and thus represent the same number. But to get the construction to work, we need as many large sets as possible. For every set of indices, either it, or its complement must be designated large. And we must do this in a consistent way that ensures every set containing a large set is large, and the intersection of two large sets is large.

A set of subsets of $\mathbb{N}$ that meets these conditions is called a *free ultrafilter*. There are many ways to add subsets so as to consistently meet these rules, but the axiom of choice is required to construct a particular example. By specifying the large sets of indices, the free ultrafilter, $\mathcal{F}$, tells us how to partition the space of sequences of real numbers into equivalence classes, each of which will act as an individual hyperreal. In some sense, every free ultrafilter defines a different version of the hyperreals. But fortunately, if the continuum hypothesis holds then the properties of the hyperreals are independent (up to isomorphism) of the particular ultrafilter used.

For each real number $r$, there is a corresponding hyperreal $^*r$ defined as the equivalence class containing $\langle r, r, r, \ldots\rangle$ (denoted $[\langle r, r, r, \ldots\rangle]$). But there are also more exotic hyperreals such as $[\langle 1, 2, 3, \ldots\rangle]$ — an infinite number denoted $\omega$ which will play a key role in the theory of hyperreal summation and integration.[1] There are also infinitesimal numbers such as $[\langle 1, {}^1/_2, {}^1/_3, \ldots\rangle]$, which is equal to ${}^1/_\omega$. And there is a great diversity of other numbers such as:

$$-\omega^2, \quad -\omega+5, \quad -\frac{1}{2\omega}, \quad \frac{1}{\omega^2}, \quad \pi-\frac{1}{\omega}, \quad \sqrt{\omega}, \quad \frac{\omega}{2}, \quad \omega+\frac{1}{\omega}, \quad \pi\omega, \quad \omega^3, \quad e^\omega, \quad \omega!$$

For any function $f$ that can be applied to a real number, we can create a transferred version, $^*f$, that applies to any hyperreal, $x$ (= $[\langle x_1, x_2, x_3, \ldots\rangle]$), by defining $^*f(x) = [\langle f(x_1), f(x_2), f(x_3), \ldots\rangle]$. The same is true for functions of more than one argument and for properties and relations. When there is no danger of ambiguity, we usually drop the $^*$ symbol. For example, the hyperreal version of $+$ is defined as $x + y = [\langle x_1 + y_1, x_2 + y_2, x_3 + y_3, \ldots\rangle]$.

A key property of the hyperreals is that they satisfy the *transfer principle*. This says that every first-order statement which is true for the reals is true for the hyperreals too. For example, because + is commutative on the reals, its hyperreal version is commutative on the hyperreals.

[1] While it is named by analogy to the $\omega$ of the ordinal numbers, these are not the same thing.

## Hyperreal summation and integration

To create a hyperreal theory of summation and integration, we simply generalise the finite sum and definite integral to allow their bounds ($a$ and $b$) to take hyperreal values:[2]

$$*\sum_{i=a}^{b} f(x) \stackrel{\text{def}}{=} \left[\left\langle \sum_{i=a_1}^{b_1} f(x), \sum_{i=a_2}^{b_2} f(x), \sum_{i=a_3}^{b_3} f(x), \ldots \right\rangle\right]$$

$$*\int_{a}^{b} f(x)\ dx \stackrel{\text{def}}{=} \left[\left\langle \int_{a_1}^{b_1} f(x)\ dx, \int_{a_2}^{b_2} f(x)\ dx, \int_{a_3}^{b_3} f(x)\ dx, \ldots \right\rangle\right]$$

Then we redefine how sums and integrals up to infinity are to be understood, replacing the standard limit-based definitions with these new constructs applied to the hyperreal number $\omega$:

$$\sum_{i=a}^{\infty} f(x) \stackrel{\text{def}}{=} *\sum_{i=a}^{\omega} f(x)$$

$$\int_{a}^{\infty} f(x)\ dx \stackrel{\text{def}}{=} *\int_{a}^{\omega} f(x)\ dx$$

We can similarly redefine integrals that go to $-\infty$ or end at a singular point, $s$:

$$\int_{-\infty}^{b} f(x)\ dx \stackrel{\text{def}}{=} *\int_{a}^{\omega} f(x)\ dx$$

$$\int_{a}^{s} f(x)\ dx \stackrel{\text{def}}{=} *\int_{a}^{s-\frac{1}{\omega}} f(x)\ dx$$

$$\int_{s}^{b} f(x)\ dx \stackrel{\text{def}}{=} *\int_{s+\frac{1}{\omega}}^{b} f(x)\ dx$$

That's it.

This new theory of summation and integration produces extremely intuitive answers to familiar sums and integrals:

$$1+1+1+\ldots = \sum_{i=1}^{\infty} 1 = *\sum_{i=1}^{\omega} 1 = \left[\left\langle \sum_{i=1}^{1} 1, \sum_{i=1}^{2} 1, \sum_{i=1}^{3} 1, \ldots \right\rangle\right] = [\langle 1,2,3,\ldots\rangle] = \boldsymbol{\omega}$$

$$2+2+2+\ldots = \sum_{i=1}^{\infty} 2 = *\sum_{i=1}^{\omega} 2 = \left[\left\langle \sum_{i=1}^{1} 2, \sum_{i=1}^{2} 2, \sum_{i=1}^{3} 2, \ldots \right\rangle\right] = [\langle 2,4,6,\ldots\rangle] = \mathbf{2}\boldsymbol{\omega}$$

$$1+2+3+\ldots = \sum_{i=1}^{\infty} i = *\sum_{i=1}^{\omega} i = \left[\left\langle \sum_{i=1}^{1} i, \sum_{i=1}^{2} i, \sum_{i=1}^{3} i, \ldots \right\rangle\right] = [\langle 1,3,6,\ldots\rangle] = \frac{\boldsymbol{\omega}(\boldsymbol{\omega}+\mathbf{1})}{\mathbf{2}}$$

---

[2] This generalisation of the finite sum is common and usually called the *hyperfinite sum*. By analogy we could call the generalised integral the *hyperdefinite integral*.

Note how the value of $2 + 2 + 2 + \ldots$ is exactly twice the size of $1 + 1 + 1 + \ldots$ and how the value of $1 + 2 + 3 + \ldots$ is the triangular number formula applied to $\omega$. Because the starred sum $({}^*\Sigma)$ is defined element-wise in terms of the finite sum, it inherits all its first-order properties such as linearity and additivity. And whenever a first-order property (like the value of a sum satisfying the triangular number formula) holds for finite sums up to every finite bound, it will hold at infinity too.

Integrals behave similarly:

$$= \int_0^\infty 1\ dx = *\int_0^\omega 1\ dx = \left[\left\langle \int_0^1 1\ dx, \int_0^2 1\ dx, \int_0^3 1\ dx, \ldots \right\rangle\right] = \boldsymbol{\omega}$$

$$= \int_0^\infty 2\ dx = *\int_0^\omega 2\ dx = \left[\left\langle \int_0^1 2\ dx, \int_0^2 2\ dx, \int_0^3 2\ dx, \ldots \right\rangle\right] = \boldsymbol{2\omega}$$

$$= \int_0^\infty x\ dx = *\int_0^\omega x\ dx = \left[\left\langle \int_0^1 x\ dx, \int_0^2 x\ dx, \int_0^3 x\ dx, \ldots \right\rangle\right] = \boldsymbol{\frac{\omega^2}{2}}$$

They also inherit first-order properties such as linearity and additivity. Perhaps most importantly, they inherit the Fundamental Theorem of Calculus.

**Theorem.** The Fundamental Theorem of Calculus applies to the starred integral.

**Proof:**

$$\begin{aligned}
*\int_a^b f(x)\ dx &= \left[\left\langle \int_{a_1}^{b_1} f(x)\ dx, \int_{a_2}^{b_2} f(x)\ dx, \int_{a_3}^{b_3} f(x)\ dx, \ldots \right\rangle\right] \\
&= [\langle F(b_1) - F(a_1), F(b_2) - F(a_2), F(b_3) - F(a_3), \ldots \rangle] \\
&= [\langle F(b_n) \rangle] - [\langle F(a_n) \rangle] \\
&= F(b) - F(a)
\end{aligned}$$

**Corollaries:**

$$\int_a^\infty f(x)\ dx = F(\omega) - F(a)$$

$$\int_{-\infty}^b f(x)\ dx = F(b) - F(-\omega)$$

$$\int_a^s f(x)\ dx = F(s - \tfrac{1}{\omega}) - F(a)$$

$$\int_s^b f(x)\ dx = F(b) - F(s + \tfrac{1}{\omega})$$

We can use the Fundamental Theorem of Calculus to solve improper integrals:[3]

$$= \int_1^{\infty} \frac{1}{x}\, dx = {*}\!\int_1^{\omega} \frac{1}{x}\, dx = F(\omega) - F(1) = \log(\omega) - 0 = \boldsymbol{\log(\omega)}$$

$$= \int_0^{1} \frac{1}{x}\, dx = {*}\!\int_{\frac{1}{\omega}}^{1} \frac{1}{x}\, dx = F(1) - F\left(\tfrac{1}{\omega}\right) = \log(1) - \log\left(\tfrac{1}{\omega}\right) = \boldsymbol{\log(\omega)}$$

We can even directly solve integrals where both bounds are improper:

$$= \int_0^{\infty} \frac{1}{x}\, dx = {*}\!\int_{\frac{1}{\omega}}^{\omega} \frac{1}{x}\, dx = F(\omega) - F\left(\tfrac{1}{\omega}\right) = \log(\omega) - \log\left(\tfrac{1}{\omega}\right) = \boldsymbol{2\log(\omega)}$$

This new theory of summation and integration *refines* the standard theory:[4]

1. Whenever the standard theory assigns a real number $r$ to a sum or integral, the hyperreal theory assigns it $r$ or a value infinitesimally close to $r$.
2. When the standard theory assigns a sum or integral the value $+\infty$, the hyperreal theory may assign it a more precise positive infinite value (and respectively for negative values).

As (1) suggests, convergent sums and integrals are often assigned values that are infinitesimally different to their standard values. For example:

$$\sum_{i=1}^{\infty} 2^{-i} \;=\; 1 - \frac{1}{2^{\omega}}$$

Whether this infinitesimal discrepancy is a bug or feature depends on the application. The usual proofs that this sum equals 1 rely on the principle that the sum must be within $\epsilon$ of 1 for every real $\epsilon > 0$ — i.e. infinitesimally close to 1. In the real numbers the only such number is 1 itself. But in the hyperreals, it is perfecty consistent for this to sum to a number infinitesimally close to 1. We might see this as generalising the idea that summing finitely many terms leaves a finite discrepancy from 1 — summing infinitely many terms leaves a corresponding infinitesimal discrepancy.

[3] Throughout the paper log without a subcript means the natural logarithm.

[4] In the jargon of alternative summation methods, hyperreal summation is *linear*, *finitely reindexable*, almost *regular* (disagreeing by at most an infinitesimal), but not *stable*. Given a particular ultrafilter, it is maximally *strong* (summing all sequences).

In some applications this discrepancy is meaningful; indeed essential. For example, the value above can be used to represent the probability that an infinite sequence of coin tosses comes up tails at some point — a near certainty, but not guaranteed. Along these lines, a promising non-standard probability theory that includes infinitesimal probabilities can be developed.

In contexts where the infinitesimal discrepancies are unwanted, they can easily be rounded off. There is a unique closest real to every finite hyperreal (its *shadow*) and we can optionally adjust the definition of the hyperreal sum and integral to take the shadow.[5] This version of the theory hews even closer to the standard theory.

In cases where the standard sum or integral diverges, the precise value of the hyperreal sum or integral sometimes depends on which ultrafilter was used to construct the hyperreals. For example, consider the poorly behaved sum $1 - 2 + 3 - 4 + \ldots$. On the standard theory this diverges, but not to $+\infty$ nor to $-\infty$. Since the partial sums alternate in value but increase in absolute size, the standard sum of the series is undefined. What does the hyperreal theory say?

$$1 - 2 + 3 - 4 + \ldots \;=\; \sum_{i=1}^{\infty} i(-1)^{i-1} \;=\; [\langle 1, -1, 2, -2, 3, -3, \ldots \rangle]$$

Which hyperreal is this? If the set of odd numbers is in the ultrafilter, then this is the same number as $[\langle 1, 1.5, 2, 2.5, 3, 3.5, \ldots \rangle] = \omega/2 + 1/2$ (since both representative sequences are identical on their odd elements). If the set of odd numbers isn't in the ultrafilter, then the evens must be. In that case, that the sum is equal to $[\langle -0.5, -1, -1.5, -2, -2.5, -3, \ldots \rangle] = -\omega/2$.

Less badly behaved sums and integrals can also have results that depend on the ultrafilter. For example, the integral:

$$\int_0^{\infty} 1 + \sin(x)\ dx$$

standardly diverges to $+\infty$. On the hyperreal theory of integration, it can take on values ranging from $\omega$ to $\omega + 2$, depending on the ultrafilter. In both cases, the ultrafilter-dependence is being caused by a kind of oscillatory behaviour in the partial sums.

There are several ways to proceed. The most modest option is to say that any sum or integral whose value depends on the ultrafilter is undefined (or, if it is determinately infinite, that it equals $+\infty$ or $-\infty$). Then there is no more problem than on the standard theories, which also decline to say any more than this about such sums and integrals.

The boldest option is to use a particular ultrafilter. Each ultrafilter gives its own theory of summation and integration with its own consistent set of answers to such sums — which one we use is up to us. This is consistent, but not entirely satisfactory. Since the ultrafilter cannot be finitely specified, we can't quite choose a specific theory and are left with some arbitrariness and unknowability in the values of sums and integrals. But given an ultrafilter, we do get a theory on which *every* sum and integral evaluates to some particular hyperreal.

[5] This can be done for infinite hyperreals too, defining shadow of $x$ (denoted $\underline{x}$) in terms of standard part (which only applies to finite hyperreals): $\underline{x} = \mathrm{st}(x - \lfloor x \rfloor) + \lfloor x \rfloor$.

Between these modest and bold options is a middle path where we divide claims about infinite summation and integration into those that are determinately true (true for all ultrafilters), determinately false (false for all ultrafilters), or indeterminate (true for some; false for others). For example, if we denote the sum $1 - 2 + 3 - 4 + \ldots$ as $S$, it is determinately true that:

- $S$ is not finite
- $|S| < \omega$
- $S \in \{ -\omega/2 , \omega/2 + 1/2 \}$
- $S - S = 0$

But it is indeterminate whether $S > 0$.

This is the approach I'll adopt in this paper.

Note that even in cases like these where the value of a sum (or integral) is ultrafilter-dependant, the change in value coming from modifications of that sum usually doesn't depend on the ultrafilter. For example, if we add 7 to one of the terms in the sum $S$, calling this new sum $S'$, then $S' > S$ on all ultrafilters and $S' - S = 7$ on all ultrafilters.

It may be possible to go further and resolve some (or all?) of the indeterminacies in a non-arbitrary way. For example, intuitively half the ultrafilters would include the odd numbers and half the even numbers. If so, we might be able to form an average across all ways of specifying the hyperreals. In this case, the average value for our alternating sum $S$ is $1/4$. Interestingly, this is precisely the value assigned to this series by Abel summation, Euler summation and Borel summation. The same holds for several other divergent series and suggests there might be a way of simultaneously generalising those theories of summation, assigning fine-grained infinite values where appropriate, and avoiding some of the indeterminacy. However, there are challenges in defining an appropriate average, and I will have to leave this as a suggestive possibility.

Before we proceed with applications of the theory, it is worth noting that this is an unusual use of the hyperreals. While they were developed to put Leibniz's infinitesimal calculus on a rigorous footing and are frequently applied to summation and integration, the field of non-standard analysis doesn't assign hyperreal values to sums or integrals. Instead, non-standard analysis uses hyperreals as an alternative means for proving the standard results of analysis (which are expressed in terms of extended reals). Because it will ultimately round off finite results to the nearest real number and infinite results to $+\infty$ or $-\infty$, non-standard analysis avoids the need to choose a particular hyperreal (such as $\omega$) for the bounds on its sums and integrals and avoids the ultrafilter-dependence. However, it also passes up the opportunity to go beyond the standard results and to *refine* or *extend* the theory of analysis.

## Numerosity

The hyperreal theory of infinite summation also provides us with a useful way of measuring the sizes of infinite sets of numbers. We can assign each set of integers a hyperreal value, called its *numerosity* (Benci & Di Nasso 2019, Błaszczyk 2021), by simply summing the indicator function for the set. For example, we have seen that the numerosity of the positive integers is:

$$\sum_{i=1}^{\infty} 1 = \omega$$

The numerosity of the entire set of integers is:

$$\sum_{i=-\infty}^{\infty} 1 = \left[\left\langle \sum_{i=-1}^{1} 1, \sum_{i=-2}^{2} 1, \sum_{i=-3}^{3} 1, \dots \right\rangle\right] = [\langle 3, 5, 7, \dots \rangle] \ = \ 2\omega + 1$$

Which is equal to the numerosity of the positive integers ($\omega$), plus the numerosity of the negative integers ($\omega$), plus the numerosity of the set containing zero (1).

More generally, to find the numerosity of a set of integers $A$, we take the hyperreal sum of the indicator function for $A$:

$$\sum_{i=-\infty}^{\infty} \chi_A(i), \quad \text{where } \chi_A(i) = \begin{cases} 1 & ,i \in A \\ 0 & ,i \notin A \end{cases}$$

The numerosity of the set of even positive integers is thus $[\langle 0, 1, 1, 2, 2, \dots \rangle] = \lfloor \omega/2 \rfloor$ while the numerosity of the set of odd positive integers is $[\langle 1, 1, 2, 2, 3, \dots \rangle] = \lceil \omega/2 \rceil$.[6] The numerosity of the square numbers is $\lfloor \sqrt{\omega} \rfloor$ while the numerosity of the non-square numbers is $\lceil \omega - \sqrt{\omega} \rceil$. More generally, if $A$ is a strict subset of $B$, then $A$ has a lower numerosity than $B$. And if a collection of subsets $A_i$ partition a set $B$, then the numerosity of $B$ equals the sum of the numerosities of the $A_i$.

These are intuitively desirable properties for a conception of size, but we have become accustomed to their absence. The familiar way of measuring set sizes via bijections, assigns all these sets the same cardinal number, $\aleph_0$. Neither way is more correct — they just represent different aspects of our intuitive conception of size as applied to infinite sets. They are complementary approaches, and it is good to have both in our toolkit.

As we will see, the hyperreal theory for sums and integrals often provides different answers to the Cantorian approach of cardinal numbers. The Cantorian reliance on bijections means it deliberately ignores the internal structure of a set of numbers and so provides a more abstracted, coarse-grained, concept of size. It isn't designed to represent the intuitive idea that the set of square numbers is becoming more and more sparse — but numerosities capture this. They even tell us that the proportion of natural numbers which are square numbers is just the ratio of these numerosities, $\lfloor\sqrt{\omega}\rfloor/\omega$ — an infinitesimal fraction, yet greater than zero.

---

[6] Where $\lfloor x \rfloor$ is the floor of $x$ and $\lceil x \rceil$ the ceiling. Note that this is another minor case of indeterminacy. While it is determinate that the numerosity of the set of even numbers is $\lfloor \omega/2 \rfloor$, it is indeterminate whether this equals $\omega/2$ or $\omega/2 - 1/2$. That depends on whether $\omega$ is even or odd, which depends on the ultrafilter.

## Infinite Expectations

In probability theory and decision theory, distributions whose means diverge to infinity cause substantial problems. For example, consider the St Petersburg gamble (Bernoulli 1738), where there is a 1 in $2^n$ chance of gaining utility[7] $n$, for all positive integers $n$. Its expected value standardly diverges to $+\infty$. This same expected value is assigned to a one-half chance of receiving a St Petersburg gamble, or to a St Petersburg gamble where all utilities have been doubled, and so on. This makes it impossible to compare such divergent distributions by their expected values, and violates important axioms of decision theory (such as preferring a distribution that dominates another and preferring higher chances of an option when it is better than its alternative).

To escape this problem, many decision theorists have concluded it is irrational to allow unbounded utilities. This forces all distributions to have finite expected value, sidestepping this problem. But it does so at the high price of requiring us to change our values — even in finite cases.[8] Other decision theorists have retained unbounded utilities and developed novel systematic ways of ordering these prospects, but have given up on assigning them numerical expected values.[9]

The hyperreal theory of summation and integration can be used to assign precise infinite expectations/means to many probability distributions whose means standardly diverge to $+\infty$ or $-\infty$. One can do this by simply taking the standard formula for the mean (the sum or integral of $x$ times $f(x)$) and applying the hyperreal interpretation of summation/integration. According to the hyperreal theory, the expectation of the St Petersburg gamble is:[10]

---

[7] Decision theorists, economists, and moral philosophers have several different ways of using the term 'utility'. In this paper, I'll treat it as an interpersonally comparable *ratio-scale* quantity. I'll sometimes use 'value' as a synonym.

[8] For example, making it no longer possible to hold that creating an additional happy life with utility $u$ has the same positive value no matter how many would be created otherwise.

[9] See, for example Colyvan (2008), Easwaran (2008), Lauwers & Vallentyne (2016), and Russell & Isaacs (2020). I take my project here to be complementary to these — finding a way of assigning intuitive infinite numbers to prospects which induce an order on prospects that is consistent with those from this growing literature.

[10] We were able to do this with an infinite sum because the utility ranged over the integers. A discrete gamble that could take non-integer values would require integrating $x$ times the PMF, and a continuous gamble would require integrating $x$ times the PDF.

Note that the hyperreal method adopts a particular ordering and structure for how to do the sum or integral. It goes through each utility level in order and adds up its contribution to the EV. If you arrange them in some other way you may end up with a different, incompatible answer. For example, the St Petersburg gamble's EV is often expressed via the contribution to

$$\sum_{i=1}^{\infty} i \cdot P(i) = [\langle 0,1,1,2,2,2,2,3,\ldots \rangle] = [\langle \lfloor \log_2(n) \rfloor \rangle] \ = \ \lfloor \log_2(\omega) \rfloor$$

This is infinite, but much smaller than $\omega$. And a probability $p$ of receiving a St Petersburg gamble is assigned expected value $p\lfloor \log_2(\omega) \rfloor$ as we would hope.

The value $\lfloor \log_2(\omega) \rfloor$ is an example of what I call a *lesser infinity* — an infinite number that is less than $k\omega$ for every real-valued $k$. This is an important class of infinite value that is absent from most well-known infinite number systems. All probability distributions over finite numbers (such as the reals) have expectations that are at best lesser infinities.[11] The St Petersburg gamble (and all other distributions over the reals) are thus beaten by any positive real chance of receiving $\omega$ value. If the value at stake in Pascal's wager (say, an eternity in heaven) is taken to be $k\omega$, and its probability is $p$, then its EV is $pk\omega$. This is also an infinite value: one worth much more than the $\lfloor \log_2(\omega) \rfloor$ of the St Petersburg gamble, albeit much less than a guarantee of $k\omega$ value.

Note that this relationship between a lesser infinite expected value and a regular infinite expected value violates the continuity axiom of expected utility theory. Indeed one can construct a dense hierarchy of probability distributions with expected values of the form $\omega^k$ for $0 < k < 1$ where any non-zero real chance of a draw from a distribution with higher $k$ beats any chance of a draw from one with lower $k$. To proceed, one either needs to either give up the continuity axiom or express it in a way that allows for the probabilities to range over infinitesimal hyperreal values. Both options seem reasonable given that in this situation some alternatives really do appear to be worth infinitely more than others.

This new approach reveals a wide range of values for distributions with infinite mean. Where the standard approach treats all distributions with positive infinite mean as equivalent, on the hyperreal approach they possess an even greater range of values than do distributions with finite means. The EVs of two different prospects can have an infinite difference or even an infinite ratio. And this can be true even if the prospects just range over the reals. This suggests that treating all prospects with infinite mean as equivalent could be as much of a mistake as treating all prospects with finite mean as equivalent.

One can also use hyperreal summation and integration to assign fine-grained infinite values to higher moments of a distribution, such as its variance. For example, the Cauchy distribution (whose mean and variance are standardly undefined) can be shown to have mean zero and variance $\approx 2\omega/\pi - 1$. This may allow one to analyse such distributions using statistical tools that standardly break when the variance is infinite. This could be important since the difficulty in working with some of these distributions often leads to them being replaced for reasons of convenience with

---

EV of $n$ successful coin flips, giving $2/2 + 4/4 + 8/8 + \ldots = 1 + 1 + 1 + \ldots = \omega$. Whereas summing the contribution of each utility level gives $0/1 + 2/2 + 0/3 + 4/4 + 0/5 + 0/6 + 0/7 + 8/8 + \ldots = 0 + 1 + 0 + 1 + 0 + 0 + 0 + 1 + \ldots = \lfloor \log_2(\omega) \rfloor$.

[11] There is thus a connection between lesser infinities and the idea of incomplete infinities: the concept of being unlimited or unbounded, yet not properly infinite.

distributions that have thinner tails — and these systematically underestimate outlier events.[12]

## Infinite Futures

In economics it is common to evaluate policies via the stream of utilities they would produce at future times. For finite streams a natural (and temporally impartial) method is to take the sum of these utilities (or in continuous time, the integral). But for unbounded streams (represented as infinite sequences of real numbers) this typically diverges to infinity, causing problems for the analysis. For example the streams of utilities given by $\langle 1, 1, 1, \ldots \rangle$ and $\langle 2, 2, 2, \ldots \rangle$ both standardly sum to $+\infty$ and so would be considered equal despite the latter being better at all times.

This problem of infinities is one of the standard motivations for exponential discounting of utility itself, also known as *pure time preference* (Koopmans 1960). For if we assume a finite bound on the level of utility at each time, then exponential discounting guarantees a finite bound on the infinite sum of the stream of utilities. The general study of how to evaluate or compare such unboundedly long streams of utility is known in economics as *intergenerational equity* (Diamond 1965, Asheim 2010, Pivato and Fleurbaey 2024).

While the infinities here may seem to stem from an implausible modelling assumption (a literally infinite future), they tend to reappear even on more plausible models. For instance, they reappear unless one is absolutely certain the future is finite. Even if one is certain that the future is finite, if its duration is probabilistic, the a declining hazard rate can make the *expected* duration infinite. And even if we are certain that there is a constant hazard rate (such that the objective expected duration is finite), uncertainty about how low this rate is can easily make the subjective expected duration infinite once more.[13] So it is not unreasonable for economists to see if one can address the infinite utility stream at the heart of the problem.

By simply taking the hyperreal sum of the sequence of utilities, we can evaluate infinite streams of utility without discounting.[14] This assigns the value $\omega$ to the stream $\langle 1, 1, 1, \ldots\rangle$ and $2\omega$ to $\langle 2, 2, 2, \ldots\rangle$. More generally, if all utilities in a stream are multiplied by $m$, the value of the stream is multiplied by $m$, and if $c$ is added to every utility, the value of the stream is increased by $c\omega$. If $c$ is instead added to a single utility, the value of the stream is increased by $c$. For instance, the stream $\langle 0, 1, 1, \ldots\rangle$ is assigned a value of $\omega - 1$. The method is thus sensitive to changes to any element of the stream and satisfies the following:

> ***Strong Pareto***:
> For all streams $x$, $y$: if $x_n \geq y_n$ for all $n \in \mathbb{N}$, and $x_n > y_n$ for some $n \in \mathbb{N}$, then $x$ is better than $y$.

---

[12] The lack of finite means and variance doesn't make such distributions less realistic — all the values they can take are still finite and most widely-used used distributions also have unbounded ranges. But a feature of these distributions is that they are very sensitive to what happens in the distant tails (e.g. to whether the upper limit is really infinity, or perhaps just $10^{100}$).

[13] This occurs when the uncertainty is uniform over some interval of the form $(0, K]$ as well as for many other distributions that aren't bounded away from zero. The uncertainty about the hazard rate produces an effective discount rate that declines over time towards the infimum of the plausible hazard rates, for reasons similar to those discussed by Weitzman (1998) (explicated in Ord (2023)). The case of uniform uncertainty for a constant hazard rate in $(0, K]$ gives an expected duration on the order of $\log(\omega)$ years.

[14] This has been proposed and rigorously explored by Pivato in a 2008 working paper. I independently suggested the same approach (in the context of infinite ethics) to Bostrom in 2004, which eventually appeared in Bostrom (2011).

At the same time it satisfies a form of temporal impartiality:

> ***Finite Anonymity***:
> Reordering finitely many elements of a stream won't change its value.

It also satisfies other desirable properties from the literature on intergenerational equity: there is no dictatorship of the present or of the future, it can handle exponentially growing utility streams, and it exhibits only infinitesimal impatience (Pivato 2008).

It is instructive to compare the hyperreal method to that of taking the standard sum and taking the standard average. The standard sum works well until the total is infinite, when it is no longer able to distinguish between infinitely valuable streams. It is like a camera whose sensor is washed out with the brightness of the scene, clipping the values to pure white. It is possible to take the average instead (Pivato 2022), allowing us to clearly see the infinite difference in value between $\langle 1, 1, 1, \ldots\rangle$ and $\langle 2, 2, 2, \ldots\rangle$. But then we can no longer discern the finite difference between $\langle 1, 1, 1, \ldots\rangle$ and $\langle 0, 1, 1, \ldots\rangle$. Taking the average is like putting a dark filter over the camera's lens — it overcomes the dazzling brightness but limits our ability to make fine distinctions in the shadows.

We say cameras have a limited *dynamic range* — the ratio between the largest and smallest values they can measure. The dynamic range of the standard sum is unbounded yet finite. It can express arbitrarily large finite ratios, but it is unable to distinguish between infinite ratios. We can adjust where we are using that dynamic range (e.g. by taking an average), but we can't increase the dynamic range without moving to a more expressive number system.

The hyperreals have a much larger dynamic range. They can handle infinite values while still keeping track of finite differences. They can even handle higher infinities, such as assigning $2^{\omega} - 1$ to the stream $\langle 1, 2, 4, \ldots\rangle$, while still being able to track finite (or even infinitesimal) differences. It is also possible to express these values as averages if we divide them by the count of the sequence ($\omega$). So the average value of $\langle 1, 1, 1, \ldots\rangle$ is 1 and the average value of $\langle 0, 1, 1, \ldots\rangle$ is $^{\omega-1}/_{\omega}$, which equals $1 - {}^{1}/_{\omega}$, and is thus infinitesimally lower. The presence of infinitesimals preserves the infinite dynamic range. Indeed, one can see this hyperreal average as a refinement of the standard average, in the same way that the hyperreal sum and integral are refinements of the standard sums and integrals.

Hyperreal summation also allows us to distinguish streams whose value, while infinite, is a lesser infinity. For example, the stream ⟨1, 0, 0, 1, …⟩ (where only times that are square numbers have utility 1) is becoming ever-more sparse as time goes on. It is assigned the value $\lfloor\sqrt{\omega}\rfloor$. There is a dense set of such lesser infinites between 0 and $\omega$, and the hyperreal method is able to distinguish all of them.[15]

The ordering on utility streams induced by hyperreal summation is closely related to the partial order induced by the well-known *overtaking criterion*. This criterion doesn't assign values to streams, but provides a method for comparing two streams $x$ and $y$ . It says $x$ is better than $y$ whenever there is a time after which the partial sum of the initial segment of $x$ remains above that of $y$:

$$x \succ y \ \text{ iff } \ \exists T \ \forall t > T : \ \sum_{i=1}^{t} x_i > \sum_{i=1}^{t} y_i$$

If we are only saying that one stream is (determinately) better than another when its hyperreal sum is higher on all ultrafilters, then the hyperreal approach's partial ordering of utility streams is exactly that of the overtaking criterion. If we instead evaluate streams according to a particular ultrafilter, $\mathcal{F}$, then the complete ordering induced will be a self-consistent refinement of the partial order given by the overtaking criterion.

Thus if all you are interested in is the ordering of utility streams, it is probably simpler to use the overtaking criterion. But if you wish to consider the differences or ratios of their values, or to compare uncertain prospects, or even just to know how valuable a utility stream actually is, you need to go beyond the overtaking criterion.

These two different interpretations of hyperreal summation lead to orderings of utility streams that are either *incomplete* (failing to rank some pairs) or *ineffable* and *arbitrary* (due to the reliance on the axiom of choice to complete the ordering). It would be good to suffer neither of these limitations, but work by Zame (2007), Lauwers (2010), and Dubey (2011) has shown that is not possible. All rankings of utility streams which satisfy Strong Pareto and Finite Anonymity are either incomplete or require the axiom of choice.

Unlike the finite sum, the hyperreal sum can be changed by adding zeros. For example, consider the stream ⟨1, 2, 3, …⟩. If we insert a neutral initial time period, moving everything else one time period along to make room, we get the stream ⟨0, 1, 2, 3, …⟩. But the hyperreal sum of this stream is $\omega$ units lower than the first stream. Shifting everything forward in time has made things worse — violating some interpretations of temporal impartiality. However, this is a direct consequence of satisfying any version of Pareto. For ⟨0, 1, 2, 3, …⟩ is also the stream we get if we subtract 1 unit of utility from every period, which Pareto requires us to rank as worse.[16] There isn't enough structure to the situation to enable us to treat these cases differently, so treating the shifted sequence as worse is a necessary cost to developing a theory that satisfies Pareto.

[15] One can create simple systems with infinite dynamic range using tuples of real numbers, where the latter places lexically dominate the earlier ones (e.g. using ($a$, $b$, $c$) where the hyperreal method would have $a + b\omega + c\omega^2$). However, such systems can't represent the dense set of lesser infinities.

While finite permutations cannot change the value of a utility stream, infinite permutations can. But again, it is not possible to satisfy Pareto while being invariant to infinite permutations (Lauwers 1998). One can't have it all, so the hyperreal approach makes a creditable attempt at taking the Pareto path and seeing what kind of a theory can be established there.

## Infinite Populations

In ethics, attempts to evaluate outcomes in terms of their total utility run into divergent sums when considering infinite worlds, containing infinite populations (Cain 1995, Vallentyne & Kagan 1997, Bostrom 2011, Askell 2018). Some philosophers are concerned about this even if our universe is finite (since moral theories should apply in all situations, even hypothetical ones) but it draws additional impetus from the fact that some mainstream cosmological theories suggest a spatially infinite universe containing infinitely many inhabited planets and thus infinitely many morally relevant beings (Bostrom 2011, Tegmark 2014).

This can be conceptualised in terms of an infinite set of persons with their own utilities, or as an infinite 3D space throughout which utility is distributed. I can't see how to apply the hyperreal method to the former, so will restrict the analysis to cases where the utility is distributed through space and a standard integral over that space diverges to infinity.

Let's start by looking at how we could solve this if there were a finite amount of utility distributed over the infinite space. We could find the amount of utility by setting up a surface and sweeping it through the space, integrating all the utility that it encounters. We could set up a symmetrical reference volume around the origin (such as a cube or sphere) and grow its size, until it surface sweeps through every point in space. If $r$ is the radius of this volume and $\rho(r)$ is the average density of utility on its surface, then we can integrate these expanding shells as they sweep through all of 3D space:

$$\int_0^\infty \rho(r) \cdot 24r^2 \; dr \qquad\qquad \int_0^\infty \rho(r) \cdot 4\pi r^2 \; dr$$

To apply the hyperreal theory, we just take the hyperreal values for these integrals.

For example, if the density of utility were constant these would integrate to $8\rho\omega^3$ and $^4/_3\ \pi\rho\omega^3$ respectively. Those constants are meaningful: space divides into 8 octants, each with value $\rho\omega^3$ for a total of $8\rho\omega^3$, while $^4/_3\ \pi\omega^3$ is the volume of a sphere of radius $\omega$. But we could also normalise these results: dividing by the size of the reference volume of unit radius (8 and $^4/_3\,\pi$ respectively) giving a normalised value of $\rho\omega^3$ in both cases.

[16] Even when the total utility over time is bounded (e.g. ⟨1, $^1/_2$, $^1/_4$, …⟩) delaying a stream by a time period and adding an initial zero can make the hyperreal sum smaller, albeit only infinitesimally so. One could see this as a kind of infinitesimal impatience.

To avoid complexities of expressing utility as a density, we could also directly construct this integrated utility from its component definite integrals as $[\langle v(R_1), v(R_2), v(R_3) \ldots\rangle]$, where $v(R_n)$ is the utility contained within the region of radius $n$.

As with streams of utility, we get the results that multiplying all utilities by $m$ multiplies the total by $m$ and adding $c$ units of utility to some location increases the total by $c$. This time, adding $c$ units to the density of utility at every location adds $c\omega^3$ units to the normalised total. We can also calculate the average utility across all space. If the utility in any region of unit volume is bounded, then the average will be a finite number which is nudged up or down by infinitesimal amounts when we make finite changes to the distribution of utility. Once again, we have the infinite dynamic range needed to keep track of both infinite and finite changes.

Unfortunately, these results require a privileged origin from which to begin our integration. This is more of an issue for space than it was for time as it is widely held that the physical space of our universe has no privileged origin. Indeed, this is just one of a host of issues for evaluating infinite spaces that arise from the complex nature of our spacetime as described by special and general relativity. So the hyperreal theory can open the door to assessing infinite worlds, but would need substantial further development before it could work in a physically realistic universe.

## Challenges & limitations

While the hyperreal theory of summation and integration shows substantial promise, this current version that I've articulated faces a number of challenges and limitations:

1. *Ultrafilter dependence*: The values of some hyperreal sums and integrals depend on the choice of ultrafilter. Possible solutions include treating such cases as undefined, choosing a particular ultrafilter, distinguishing between claims that are true on all/some/no ultrafilters, or finding a way of averaging across ultrafilters.

2. *Grid problem*: The hyperreal value of an improper integral only depends on the values of definite integrals at a countable sequence of points (either the integers or on numbers of the form $s \pm {}^1/_n$, where $s$ is a singular point). This ignores any divergent behaviour occurring between these points. Possible responses include characterising the functions where this is an issue and restricting the scope of this method to not apply to them, replacing the hyperreals with a system that doesn't have this problem, or just accepting this dependence.

3. *Why ω?*: We currently lack a principled reason for choosing the hyperreal $\omega$ as the unit by which to measure the infinite. Possible responses include finding such a reason, using a different number system (like the surreals whose $\omega$ *is* more special), or arguing that it doesn't matter which hyperreal we use and thus choosing $\omega$ (or any other number) by fiat to act as the unit of the infinite (much as we could have chosen any platinum rod to act as the metre).[17]

4. *Why 1/ω?*: Is there a principled reason for using the reciprocal of the infinite bound for determining how close to a singular point a hyperreal bound should come? One alternative would be to deal with vertical asymptotes by integrating against the $y$-axis instead of the $x$-axis, though this gives different values to the integrals, and loses some of the properties we've established. Another alternative would be to integrate up to the $x$-coordinate where the $y$-coordinate reaches $\omega$ (i.e. to $s \pm f^{-1}(\omega)$ instead of $s \pm {}^1/_{\omega}$)).

5. *Spillover*: Sums and integrals to infinity are supposed to be summing or integrating over all the positive integers or reals. Summing or integrating up to an infinite hyperreal bound like $\omega$ ensures we cover them all. But stopping at ${}^{\omega}/_2$ would also achieve this. Including contributions to the total from infinite values of $x$ between ${}^{\omega}/_2$ and $\omega$ seems inappropriate. But since there is no least infinite hyperreal, all hyperreal stopping points would appear to include too large a domain. (An equivalent issue exists for $-\omega$ and for integrating to within ${}^1/_{\omega}$ of a singular point.)

6. *Lack of closure*: These hyperreal methods allow us to sum or integrate functions over the reals and provide a result that may not be a real number. But these methods don't directly let us sum or integrate functions of those non-real outputs. So while extending our number system to the hyperreals may solve the old problems, it also raises new problems it cannot solve. This is a particular issue when taking a sum or integral across two orthogonal dimensions, as we can't just use the familiar approach of nested sums or integrals.

A different kind of critique concerns whether the use of the hyperreals is really necessary for achieving these results. Can't we order divergent sums and integrals with other methods, and even have some quantitative measures of how large each is? Most notably, we could compare them via their asymptotic behaviour and get a conception of size via the functional form of the antiderivative, $F(x)$. Indeed, one might even ask if the hyperreal approach isn't just 'asymptotics in disguise'?

There is much I'm sympathetic to in this critique. I don't make any claim that the hyperreal approach is *necessary* for evaluating divergent sums and integrals. Indeed, if any sequence $\langle x_n \rangle$, representing $x$, overtakes a sequence $\langle y_n \rangle$, representing $y$, then $x > y$ regardless of which ultrafilter is used. So this ultrafilter-independent core of the ordering could clearly be modelled directly with asymptotics.

However, much the same could be said about real numbers. The Cauchy sequence construction represents each real as an equivalence class of sequences of rational numbers and compares them by the asymptotics of those sequences. So in some sense irrational numbers are also just asymptotics in disguise. But something substantial would be lost if we treated all claims about irrational numbers that way. e.g. if instead of saying that the area of the circle was $\pi r^2$, we had to say that despite having no (rational) area, its area asymptotically approaches that of a regular polygon of radius $r$ as the number of its sides, $n$, increases to infinity. We could make do with this, but mathematics would be needlessly harder. If we'd never treated $\pi$ as being a proper number, we couldn't use it in our equations, without wrapping every line where it appears (or could be the result) in a complicated limit.

---

[17] One might also try interpreting the integral from $a$ to $\infty$ as a starred integral from $a$ to $\omega + a$. This ensures that the integral always spans the same horizontal width ($\omega$) and makes it invariant to horizontal shifts. However, there are serious tradeoffs such as losing the property that the integral from 0 to $\infty$ equals the sum of the integrals from 0 to $a$ and from $a$ to $\infty$.

Treating $\pi$ (and the other irrationals) as bona fide numbers greatly facilitated the progress of mathematics by making mathematics more eloquent; less cumbersome. The same is true for the complex numbers. My contention is that *despite* hyperreal summation and integration being largely equivalent to asymptotics, it is nevertheless important to assign *numbers* to these sums and integrals, and that the hyperreals can perform that job.

There are some important limitations of hyperreal summation and integration as applied to decision theory, economics, and ethics.

In all these areas, it always provides the same partial ordering as the appropriate type of overtaking criterion. So it only really goes further than the overtaking criterion if you are committing to a particular ultrafilter (which refines this into a complete ordering at the cost of arbitrariness and ineffability) or if you want to *evaluate* options, not just to compare them. This could be for the purpose of using their numerical representations in further calculations, or more philosophically, to show that the paradigm of *options having values* need not be revised in light of infinite cases.

I've presented hyperreal summation and integration as an attempt to solve problems of evaluating options comprised of infinitely many parts, each with finite value. There are many other challenges related to infinities that it doesn't attempt to solve. For example, it doesn't immediately offer any way of evaluating indivisible parts of infinite value. The hyperreals may be rich enough to ascribe values to such parts, but the theory presented here doesn't show which values to assign.

Another challenge is *fanaticism* — the fact that on many theories of how to assess infinite options, even an arbitrarily small finite chance of an infinitely valuable option beats all finitely valuable options. The theories I'm outlining here are fanatical in this sense. But since one can find troubling analogues of fanaticism even in large enough finite cases (Bostrom 2009), I'm not sure the solution (if any) lies in how one deals with the infinite. My best guess is that the solution doesn't lie in the values we assign to these options, but in limits on how demanding morality can be. For example, we could allow that it would be *better* to choose a small gain for every future generation over a large gain for our own generation (and that our generation would be praiseworthy if we chose this), but that it is sometimes *permissible* to choose a suboptimal option.

There are also many challenges and paradoxes of infinite choice — problems that arise when one has infinitely many simultaneous options or an infinite succession of choices. The hyperreal theory is not put forward as an attempt to solve these, so they would need to be resolved in some other way.

In finite cases, one can easily connect the theories of how we sum across possibility, time, and space, into a coherent whole. At its simplest, we can just nest the sums or integrals, treating time, space, and possibility as different dimensions in which utility can reside (Broome 1991). But there are challenges in connecting the hyperreal sums or integrals in this way. For example, doing so would suggest that a world consisting of $\omega$ utility at each time (the utility stream $\langle\omega, \omega, \omega, \ldots\rangle$) would have a value of $\omega^2$. But this value has already been assigned to the utility stream $\langle 1, 3, 5, \ldots\rangle$ — a world that is strictly worse at all times. There is a subtle incompatibility between the idea of $\omega^2$ as the product of two infinite dimensions, and as the value of infinite quadratic growth in a single dimension. This would need to be resolved to create a fully satisfactory theory across all dimensions at once. One possible route forward is to adopt the approach we used for multiple spatial dimensions: starting with a finite 'volume' that extends in all the dimensions and to integrating it until it covers the entire space.[18]

There are also two issues that affect the hyperreal theory of infinite expectations.

First, it violates a version of the *sure-thing principle*: that the value assigned to a prospect should never be higher than that of all its possible outcomes. For example, the hyperreal method suggests a ticket for a St Petersburg gamble is worth an infinite value ($\lfloor\log_2(\omega)\rfloor$) so we should be indifferent to trading away a guaranteed outcome worth that value to get such a ticket. However, we know with certainty that after the gamble has been resolved, we will get a finite reward and regret the trade (Russell & Isaacs 2020). So any infinite valuation would be too high. But there is a corresponding problem if we assign any finite value to the gamble: for all finite valuations, there is a finite version of the St Petersburg gamble (truncated at some large number of coin tosses) that has expected value exceeding our valuation, so the full gamble must be worth more than we said.

One possible way out would be to lean on the novel fact that the hyperreal theory assigns a decidedly low infinite value to the St Petersburg gamble. We could deny that guaranteed outcomes can have lesser infinites as values, allowing us to order these divergent gambles above all finite outcomes but below all infinite ones. However, as lesser infinite values naturally come up in the study of infinite time and space, this would come at a serious cost to the larger project of a uniform approach to all of these infinities. Ultimately, my best guess is that we need to abandon this version of the sure-thing principle in infinite contexts. After all, Cain (1995) already showed that similar dominance arguments must fail in infinite cases. However, abandoning it will no-doubt have a price and is a downside of this approach.

Second, I suggested we calculate these by the hyperreal interpretation of the standard formula for the mean, where we sum (or integrate) over each utility multiplied by its probability. However, this involves summing over the domain of utility instead of summing over the domains of time or space. This is a problem because the hyperreal sum and integral are invariant to vertical translations and dilations but aren't invariant to horizontal translations and dilations. Thus, the hyperreal theory can't cleanly derive that a prospect with twice as much utility in each state is twice as good (it ends up depending on the ultrafilter).

[18] For example, if the $\omega$ utility at each time were understood as 1 unit of utility at each location on an infinite sequence of locations, then the integral formed from increasingly large space-time squares would converge to $\omega^2$, which is the same value it would reach if the utility were zero everywhere except for one location where it increased as $\langle 1, 3, 5, \ldots\rangle$. So there is a natural sense in which these should produce the same value, despite the former dominating the latter at every time.

One promising path forward is to rearrange the way we do the sum: to integrate over probability from 0 to 1, with states arranged in order of increasing utility. In other words, we take the expected value to be the hyperreal integral of the distribution's quantile function: summing up the contribution to the EV from each quantile in turn. For the St Petersburg gamble, the integral would be:

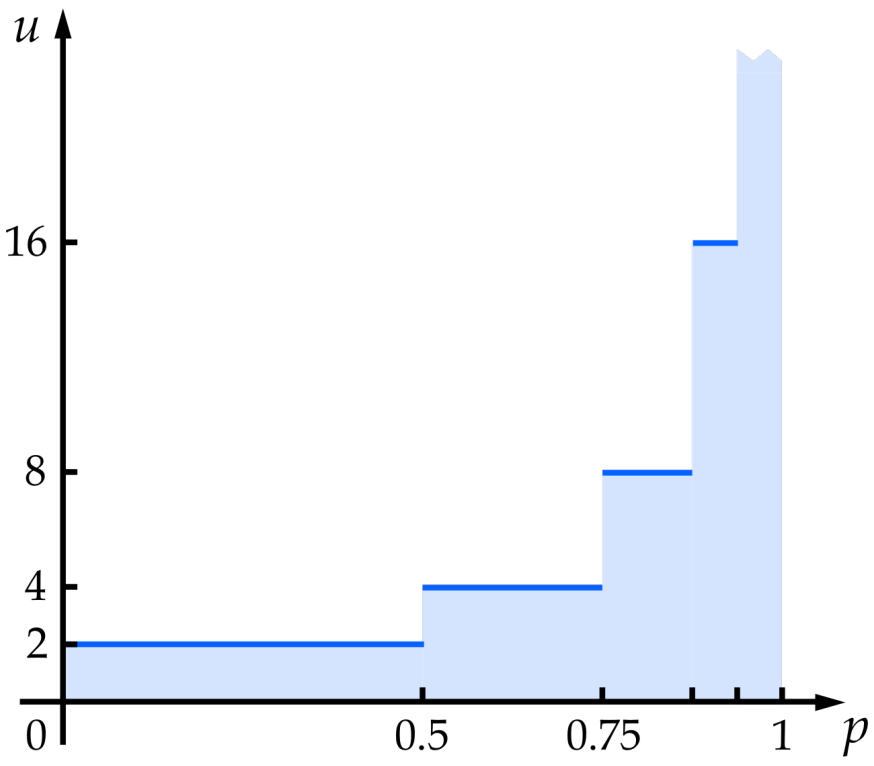

This may be able to form the basis of a hyperreal expected value that is linear in the utilities of the states as well as in their probabilities.

## Further Work

This paper has only provided a quick sketch of the mathematical theory of hyperreal summation and integration so there is much more work needed to properly flesh it out. Some natural paths for further development include:

- Characterising the set of functions whose hyperreal integrals take on values that are independent of the ultrafilter used;
- Developing the idea of averaging over ultrafilters;
- Developing versions based on other infinite number systems, such as the surreal numbers.[19]

There are also many ways in which one could develop the applications of the theory. Each of the areas we've covered could be explored much more deeply: teasing out the consequences of evaluating infinite options in this manner and comparing these with both competing theories and the known impossibility results in that area.

---

[19] Pivato (2014) explores a version based on an alternative formulation of the hyperreals.

I conjecture that the canonical answers given by the hyperreal theory will transfer to any improved theories. For instance, that the integral from 1 to $\infty$ of $^1/_x$ will equal $\log(\Phi)$ where $\Phi$ is a simple and canonical infinite number in that theory.

In each of these areas, the hyperreal approach could also be applied to other forms of aggregation beyond direct sums and integrals. For example, in ethics *prioritarianism* (Parfit 1997) gives more priority to helping the worse off by transforming people's utilities with a concave function before summing them, and *risk-weighted expected utility theory* (Buchak 2013) allows for risk aversion about utility itself through use of a more complex kind of sum across states. Many such systems can be adapted to use the hyperreal sum. Even the product of utilities can be adjusted to a hyperreal infinite product in the natural way.

And there are more areas where hyperreal summation and integration could be of use:

Quantum field theory is a very different topic that is also beset by infinities. In a variety of situations, crucial sums and integrals diverge to infinity. Physicists have developed a diverse set of techniques under the names of regularization and renormalisation to work around these issues. It is possible that hyperreal summation and integration (or some development of these theories) would help.

For example, in some notable physical situations infinities arise from the sums $1 + 1 + 1 + \ldots$ and $1 + 2 + 3 + \ldots$. Physicists have had substantial success using zeta regularisation to assign these sums the values $-{}^{1}/_{2}$ and $-{}^{1}/_{12}$, respectively. Despite the bizarreness of assigning a negative fractional value to a sum of strictly positive integers, many summation techniques assign them precisely these values.

Looking more carefully at what is happening in the physical situations, the sums appear to contain a term that is diverging to infinity (linearly for the first; sum quadratically for the second) and a term that is converging to a finite number ( $-{}^{1}/_{2}$ and $-{}^{1}/_{12}$ respectively). Other aspects of the physical situations suggest that the infinite parts cancel with another (intuitively equal) infinity, leaving the residual finite part. While this might seem fanciful, in the case of the Casimir effect the resulting $-{}^{1}/_{12}$ has been physically observed, suggesting that this really is the correct answer, even if we don't yet know how to rigorously justify it.

This looks like a situation where the infinite dynamic range of hyperreal numbers could help. For example, the situations above might be explicable if the integral related to $1 + 1 + 1 + \ldots$ took the value $\omega - {}^{1}/_{2}$ and the integral related to $1 + 2 + 3 + \ldots$ took the value ${}^{1}/_{2}\ \omega^{2} - {}^{1}/_{12}$. These hyperreal valuations do not immediately fall out of the theory presented here, but they do appear to be plausible values of the smoothed functions corresponding to the discrete infinite sums. I believe that further work here is likely to turn up a way of making this precise. This may provide some help with this longstanding technical problem of infinites in physics. (And it may also be of mathematical interest — illuminating the relationship between divergent sums and integrals.)

Hyperreal summation and integration might also be able to provide the foundation for a theory of infinitesimal probabilities. A key aim of such a theory is to avoid assigning probability zero to events that can happen (such as to each possible infinite sequence of coin flips). One way to do this is through the survival function, $S(t)$, which represents the chance that an event still hasn't occurred by time $t$. For example, let the event be that the coin has at some point landed tails. If the coin is being flipped at rate $\lambda$ starting at time $t_0$ then:

$$S(t) = \begin{cases} 1 & ,t < t_0 \\ 2^{-\lfloor r(t-t_0)+1 \rfloor} & ,t \geq t_0 \end{cases}$$

We could begin a theory of infinitesimal probability by proposing that the chance an event never happens be represented by the value of the survival function at time $\omega$ (after all finite times have passed).[20] Thus, the probability of flipping all heads (or of any other particular sequence) would be:

$$S(\omega) = 2^{-\lfloor \lambda(\omega-t_0)+1 \rfloor}$$

This theory would assign non-zero infinitesimal probabilities to all processes that have a non-zero chance of surviving for every finite duration. And these assigned probabilities have a variety of attractive properties such as:

- starting the sequence of flips one flip later doubles the probability of all heads
- doubling the rate of flips squares the probability
- the probability of heads on every odd flip is the square root of the probability, which is the same as if we halved the rate
- the probability of *all heads or all tails* is twice the probability of all heads
- the probability of rolling all 6s on an infinite sequence of die rolls is given by the same formula, but with a 6 replacing the 2
- as a continuous example, the probability of an unstable isotope never decaying is $e^{-\lambda(\omega - t_0)}$, where $\lambda$ is its decay rate
- The probability of surviving an infinite succession of hazards with probabilities $1/2$, $1/3$, $1/4$, … is $1/\lfloor\lambda(\omega - t_0) + 1\rfloor$

This theory does require a zero and unit for time, but canonical answers can be found (the Big Bang & Planck time). And for many questions about the comparative sizes of infinitesimal probabilities the particular zero and unit factor out anyway. There is clearly a lot more to develop, but a theory of infinitesimal probabilities along these lines appears promising.[21]

---

[20] This is a kind of hyperreal sum or integral, since the survival function is 1 minus the sum or integral of the probability distribution. The value of $S(\omega)$ is equal to the infinitesimal discrepancy between this hyperreal sum or integral and its standard value of 1 (i.e. to the amount of probability remaining after summing or integrating to $\omega$).

[21] One might even be able to make a whole field of non-standard statistics, transferring a host of statistical ideas to hyperreal settings. For example, take the finite uniform distribution on the interval $[a, b]$ with density $1/(b-a)$ and generalise it to hyperreal bounds. Then a uniform distribution over the interval $[0, \infty)$ might be representable as a uniform hyperreal distribution over $[0, \omega]$ with density $1/\omega$.

## Conclusions

We have explored a novel mathematical technique for summation and integration that assigns intuitive fine-grained infinite values to many divergent sums and integrals. Given the hyperreal numbers, the method is extremely simple and elegant: we just took the natural generalisation of the finite sum and integral to hyperreal bounds and then reinterpreted infinite sums and integrals as these generalised versions with a bound of $\omega$.

The resulting theory is of mathematical interest as a way of refining the concept of the infinite sum and improper integral, and it can be applied to many areas outside of mathematics where the infinities thrown up by the standard theory block our ability to analyse and compare different options. To do so, we just took the natural infinite sums and integrals that would define the relevant values and used the new method to evaluate them. This produced intuitive results that avoid many of the problems of infinity that had plagued these areas.

One upshot of this is that we should be very skeptical about applying familiar rules for the infinities of the extended reals and cardinal numbers when evaluating or comparing infinite options. For example, Hájek (2003) argued that because the expected value of Pascal's Wager is $+\infty$, the value of a 1 in a million chance of receiving the wager is also $+\infty$ and so would (surprisingly) be equally good. But while this is true if we measure infinite value with the extended reals, it is false if we measure it with the hyperreals. Since Hájek's reasoning also violates the foundational state-wise dominance principle of decision theory (Askell 2018), having a theoretical alternative that preserves this principle is welcome.

Similarly, there is a bijection between the people in an infinite universe where 99% of people on each planet are happy while 1% are miserable and the people in a universe with the opposite ratio. But this doesn't compel us to treat them as equally valuable. The *cardinalities* of happy and miserable people are the same in both cases, but it is far from clear that cardinal numbers are the appropriate infinite number system with which to evaluate this case, and the hyperreal approach allows us to reach the intuitive answer that the former universe is superior by an infinite amount.

In his groundbreaking work on the longterm prospects for humanity, Freeman Dyson (1979) suggested that an infinite future is possible, even with a finite amount of energy. By alternating fixed-length periods of activity with exponentially increasing periods of hibernation, that finite energy could produce infinitely many years of conscious experience. This is usually taken as producing a future that would be equally valuable as one that didn't have the hibernation periods. But with this new lens of hyperreal summation, we can see that even were this scheme workable, the increasingly sparse stream of utility it produces may have a value on the order of $\log(\omega)$ — quite an impoverished infinity compared to that of constant flourishing over all time.

So I suggest that we refrain from treating claims of infinite value as claims that this value obeys the mathematics of the extended reals or cardinal numbers. Instead we should take most claims of infinite value merely to be saying that the value exceeds all finite levels. How it then compares with other infinitely valuable options is a further matter, and is not immediately given by facts about the extended reals or cardinal numbers.

Indeed, the possibility of hyperreal summation and integration allows a conception of infinite value that behaves very much like finite value. Doubling a positive value makes it twice as good, doubling the chance of a positive outcome makes it twice as good in expectation, subtracting 7 units of utility from a part of an outcome makes the outcome worse (by 7 units) etc. With the high dynamic range offered by the hyperreals, you can take in the brilliance of the infinite and still see the subtle shading within it, while remaining sensitive to the variation in the darkness of the finite shadows. Much of this is due to the transfer principle, where all first-order properties of the reals are inherited by the hyperreals. Or as Leibniz (1702) put it: 'the rules of the finite are found to succeed in the infinite.' We can thus see hyperreal summation and integration as fleshing out a Leibnizian conception of infinite value, offering an alternative to the Cantorian conception.

On this view, the appropriate family of infinite numbers to describe these infinite situations (the hyperreals) is one that is less familiar, but ultimately more intuitive. The peculiar properties of the extended reals (where $x + 1$ can equal $x$) or of the cardinal numbers (as described in the story of Hilbert's hotel) come from their coarse-grained treatment of the infinite. Once we have the fine-grained conception offered by the hyperreals, a more intuitive picture of the infinite emerges.

Access to such methods should make us less afraid of the infinite — less inclined to introduce assumptions or values that would otherwise be unwarranted simply in order to make the infinities go away. For it looks like there is a way of answering these questions that avoids the decision paralysis that comes with the coarse-grained infinities. And with such prospects available, it would be rash to change our fundamental values about finite situations out of fear that they can't be made to work in infinite situations.[22] While we don't yet have a fully developed theory for evaluating infinite options using fine-grained infinite values, we have scarcely begun to explore even the most natural infinite number systems, and what little we've seen counsels hope rather than despair.

[22] Such as by requiring pure time preference or bounded utility.